\documentclass[10pt]{article}
\usepackage{amsmath,amssymb}

\newcommand{\al}{\alpha}
\newcommand{\lam}{\lambda}
\newcommand{\frb}{\mathfrak{b}}
\newcommand{\frg}{\mathfrak{g}}
\newcommand{\frsl}{\mathfrak{sl}}
\newcommand{\frso}{\mathfrak{so}}
\newcommand{\uxg}{U_\chi(\frg)}
\newcommand{\lxl}{L_\chi(\lambda)}

\begin{document}

\title{Analogues of Weyl's Formula for Reduced Enveloping Algebras}

\author{J.E. Humphreys
\\ Dept. of Mathematics \& Statistics,
University of Massachusetts,
Amherst, MA 01003
\\E-mail: jeh@math.umass.edu
}

\maketitle

\begin{abstract} In this note we study simple modules for a reduced
enveloping algebra $\uxg$ in the critical case when $\chi \in \frg^*$
is ``nilpotent''.  Some dimension formulas computed by Jantzen suggest
modified versions of Weyl's dimension formula, based on certain
reflecting hyperplanes for the affine Weyl group which might be
associated to Kazhdan--Lusztig cells.  \end{abstract}

\section{Introduction}

In the last decade or so there has been significant progress in
understanding the non-restricted representations of the Lie algebra of
a reductive group over a field of prime characteristic.  Friedlander
and Parshall extended the earlier foundations laid by Kac and
Weisfeiler, while Premet proved the Kac--Weisfeiler conjecture on the
minimum $p$-power dividing dimensions.  More recently the work of
Jantzen has reinforced ideas of Lusztig which arise in the framework
of affine Hecke algebras and Springer fibers in the flag variety.

In spite of the progress made, serious obstacles remain to a
definitive treatment of the representations.  Here we attempt to
interpret Jantzen's explicit dimension calculations in terms of
analogues of Weyl's classical formula, imitating Kazhdan--Lusztig
theory for the case of restricted representations.

\section{Reduced enveloping algebras}\label{sec.red}

First we recall briefly some essential background and notation,
referring for details to the survey \cite{hu98} and the lectures by
Jantzen \cite{ja98}, whose notation we mainly follow.  For Lusztig's
perspective on these questions, see \cite{lu01}.

\subsection{} Let $G$ be a simply connected, semisimple algebraic
group over an algebraically closed field of characteristic $p>0$, with
Lie algebra $\frg$.  Following work of Kac and Weisfeiler, the simple
modules for the universal enveloping algebra $U(\frg)$ partition into
modules for quotients $\uxg$ of $U(\frg)$ (reduced enveloping
algebras) associated with linear functionals $\chi \in \frg^*$.  All
$\chi$ in a coadjoint $G$-orbit yield isomorphic algebras.  If
$\chi=0$, $\uxg$ is the restricted enveloping algebra, whose
representations include those derived from representations of $G$.

It has been known since early work of Jacobson and Zassenhaus that the
maximum possible dimension of a simple module for $\frg$ or $U(\frg)$
is $p^N$ ($N=$ number of positive roots). The Steinberg module in the
restricted case is an example where this dimension is achieved.

\subsection{} The ``nilpotent'' $\chi$ (including $\chi=0$) play the
main role.  These correspond to nilpotent elements of $\frg$ when
$\frg$ can be identified in a $G$-equivariant way with $\frg^*$, and
form finitely many $G$-orbits (corresponding naturally to the
characteristic 0 orbits when $p$ is good).  Probably the most
important question about the representation theory of $\uxg$ is this:

\emph{Question.  How does the geometry of the $G$-orbit $G \chi$
influence the category of $\uxg$-modules?}

The orbit geometry involves a number of important ideas which have
played a major role in characteristic 0 representation theory:
Springer's resolution of the nilpotent variety, the flag variety and
Springer fibers, affine Weyl groups and Hecke algebras,
Kazhdan--Lusztig theory.  It seems clear from recent work of Lusztig
that many of these same ideas should recur in prime characteristic.
In particular, the affine Weyl group $W_p$ relative to $p$ (defined in
terms of the Langlands dual of $G$) has for a long time been known to
play a major role in organizing the representation theory of $G$.

\subsection{} The category of $\uxg$-modules can be enriched by adding a
natural action of the centralizer group $C_G(\chi)$.  When this group
contains at least a 1-dimensional torus $T_0$, Jantzen is able to
obtain graded versions of the Lie algebra actions and exploit
translation functors much as in the restricted case.

The best-behaved case occurs when $\chi$ has \emph{standard Levi form}
in the sense of Friedlander--Parshall \cite{fp90}: for some choice of
Borel subalgebra $\frb$, we have $\chi(\frb)=0$ while $\chi$ vanishes
on all negative root vectors $x_{-\al}$ except for a set $I$ of simple
roots $\al$.  (This always happens in type $A$.)  Then the simple
$\uxg$-modules are parametrized uniformly by linked weights $w \cdot
\lam$ with $w$ running over coset representatives for the subgroup
$W_I$ of the Weyl group $W$ generated by corresponding reflections.

In general the parametrization by weights is much less well
understood.  It may depend in part on the choice of a Borel subalgebra
on which $\chi$ vanishes: such a Borel subalgebra lies on one or more
irreducible components of the Springer fiber.  It is also possible
that the component group $C_G(\chi)/C_G(\chi)^\circ$ and its
characters will play a significant role, as they do in Springer
theory.  Recent work of Brown and Gordon \cite{bg01} confirms, at any
rate, that the blocks of $\uxg$ (when $\chi$ is nilpotent) are in
natural bijection with linkage classes of restricted weights.  Here we
consider only the most generic situation, involving simple modules in
blocks parametrized by $p$-regular weights.  This requires $p \geq h$
(the Coxeter number).

\subsection{} So far the most striking general fact about $\uxg$-modules
is the theorem of Premet \cite{pr95}, valid for arbitrary $\chi$
(under mild restrictions on $\frg$ and $p$):

\emph{If $d$ is half the dimension of the coadjoint orbit $G \chi$,
then the dimension of every $\uxg$-module is divisible by $p^d$.}

This had been conjectured much earlier by Kac and Weisfeiler.  In
particular, when $\chi$ is regular, all simple modules have the
maximum possible dimension $p^N$ ($N=$ number of positive roots).
Premet's theorem suggests a natural question:

\emph{With $d$ as above, does there always exist a simple
$\uxg$-module $\lxl$ of the smallest possible dimension $p^d$?}  

The answer is yes in the cases investigated so far, but for no obvious
conceptual reason unless $\chi$ lies in a Richardson orbit (permitting
an easy construction by parabolic induction from a trivial module).
Our proposed interpretation of dimension formulas stems partly from
trying to understand this question better.

\section{The restricted case}{\label{sec.res}}

We recall briefly the standard framework \cite{ja87} for the study of
simple $G$-modules, which include all simple $\uxg$-modules when
$\chi=0$.

For each dominant weight $\lam$ there is a Weyl module $V(\lam)$,
whose formal character and dimension are given by Weyl's formulas.  In
particular, \[\dim V(\lam) = \frac{\prod_{\al>0} \langle
\lam+\rho,\al^\vee \rangle} {\prod_{\al>0} \langle \rho,\al^\vee
\rangle}.\] Each Weyl module has a unique simple quotient $L(\lam)$.
Those for which $\lam$ is \emph{restricted} (the coordinates of $\lam$
relative to fundamental weights lying between 0 and $p-1$) are
precisely the $p^r$ simple $U_0(\frg)$-modules, where $r$ is the rank.
Knowing just these modules would allow one to recover all $L(\lam)$ as
twisted tensor products, by Steinberg's Tensor Product Theorem
\cite{st63}.  But so far the broader study of Weyl modules for $G$ has
yielded the most concrete results.

Knowledge of the formal characters and dimensions of all $L(\lam)$ is
equivalent to knowledge of the composition factor multiplicities of
all $V(\lam)$.  When $p<h$ (the Coxeter number), there is no specific
program for finding these multiplicities, but for $p \geq h$ the
answer is expected to be given by Lusztig's conjecture.  (This is
known to be true for ``sufficiently large'' $p$, from the work of
Andersen--Jantzen--Soergel \cite{ajs94}.)

Lusztig's approach depends on the fact that composition factors of
$V(\lam)$ must have highest weights linked to $\lam$ under the
standard dot action of the affine Weyl group $W_p$ relative to $p$.
Write dominant weights as $w \cdot \lam$, where $\lam$ lies in the
lowest alcove of the dominant Weyl chamber.  One can in principle
express the character of $L(w \cdot \lam)$ as an alternating sum (with
multiplicities) of the known Weyl characters for various weights $w'
\cdot \lam \leq w \cdot \lam$.  The multiplicities are in turn
predicted to be the values of Kazhdan--Lusztig polynomials for pairs
in $W_p$ related to $(w',w)$, after evaluation at 1.  This procedure
is inherently recursive and even in low ranks cannot usually be
expected to produce simple closed formulas for characters or
dimensions of simple modules.

Note how use of the lowest dominant alcove as a starting point locates
in a natural way the unique weight $\lam=0$ for which $L(\lam)$ has
the smallest possible dimension $p^0=1$.  This weight is as close as
possible to all hyperplanes bounding the alcove below, i.e., minimizes
the numerator of Weyl's formula.  But cancellation by the denominator
is needed to produce 1.

\section{Special cases}{\label{sec.spec}}

At the opposite extreme from $\chi=0$, in the case where the coadjoint
orbit of $\chi$ is \emph{regular} (with $d=N$), one has $\dim \lxl =
p^N$ for all $\lam$.  Much less is known between these extremes.

In a series of recent papers, Jantzen has studied a number of special
cases when $\chi$ is nilpotent and of small codimension in the
nilpotent variety.  He obtains explicit dimension formulas for simple
modules, as well as many details about projective modules, Ext groups,
etc.  Unlike the case $\chi=0$, it is feasible here to work out closed
formulas for dimensions.

\begin{itemize}

\item[(a)] Type $B_2$, with $\chi$ in the minimal (nonzero) nilpotent
orbit was first treated by \emph{ad hoc} methods in \cite{ja97} and
then more systematically in \cite{ja00}.  We take a closer look at
this in the following section.

\item[(b)] The case when $\chi$ lies in the subregular orbit ($d=N-1$)
was initially treated in \cite{ja99a} for the two simple types $A_n,
B_n$ where $\chi$ can be chosen in standard Levi form.  A more
comprehensive treatment was then given in \cite{ja99b}.  The results
are more complete for simply-laced types.  When there are two root
lengths, the number of simple modules in a typical block is less
certain (leading to uncertainty about some of the dimensions), but
everything is expected to agree with Lusztig's predictions.

\item[(c)] Unpublished work by Jantzen (assisted by B. Jessen for type
$G_2$) deals with a number of other cases, including the nilpotent
orbits of $G_2$ for which $d=3,4$ (while $N=6$) and the ``middle''
orbit of $A_3$ (with $d=4, N=6$).  He also works out families of
examples involving standard Levi form: $C_n (n \geq 3)$ with $I$ of
type $C_{n-1}$ and $D_n (n \geq 4)$ with $I$ of type $D_{n-1}$.  In
each case $d=N-2$.  The results are somewhat less complete in types
$G_2$ ($d=3$) and $C_n$, just as in the subregular case.

\end{itemize}

It is a striking fact that, in all of these cases, the dimension
formulas for simple modules have the same quotient format as Weyl's
formula.  There is a constant denominator, together with a numerator
written as the product of $N$ factors: $p$ repeated $d$ times (in
accordance with Premet's Theorem), as well as $N-d$ other factors.
Each of these factors involves an affine expression in the coordinates
of a $\rho$-shifted weight based in one reference alcove.  One or more
weights will minimize the numerator, giving a dimension equal to $p^d$
after dividing by the denominator.

The main drawback to these formulas is that there is a separate one
for each simple module (or small family of simple modules) in a
typical block.  Moreover, there is no obvious way to predict the
formulas in advance, apart from the occurrence of $p^d$.

\section{Example: $B_2$}{\label{sec.B2}}

To explain more concretely our approach to dimension
formulas, we look at type $B_2$ (say $p \geq 5$).  Denote the simple
roots by $\al_1$ (long) and $\al_2$ (short), with corresponding
fundamental weights $\varpi_1$ and $\varpi_2$.

\subsection{} Consider the case when $\chi$ lies in the minimal
nilpotent orbit (with $N=4, d=2$).  Here $\chi$ has standard Levi
form, relative to the subset $I=\{\al_1\}$. There is a one-dimensional
torus $T_0$ in $C_G(\chi)$ which acts naturally on $\uxg$-modules.  A
generic block has four simple modules $\lxl$, each labelled by two
``highest'' weights $\lam$ linked by the subgroup of $W$ generated by
the simple reflection $s_1$.  The dimensions of the $\lxl$ were first
worked out by Jantzen in \cite{ja97}; he later developed a streamlined
version based on the systematic use of translation functors in
\cite[\S 5]{ja00}.  As required by Premet's Theorem, $p^2$ divides all
dimensions.

To parametrize the simple modules by weights, Jantzen starts in the
conventional lowest alcove of the dominant Weyl chamber, fixing a
$p$-regular weight $\lam$.  In order to simplify formulas, he builds
in the $\rho$-shift by writing $\lam+\rho = r\varpi_1 + s\varpi_2 =
(r,s)$.  Thus $r,s>0$ while $2r+s<p$.  The dimensions of simple
modules corresponding to linked weights in the four restricted alcoves
are: \[\frac{s(p-2r-s)}{2} p^2, \; \frac{2pr}{2} p^2, \;
\frac{(2p-s)(p-2r-s)}{2} p^2 \; \frac{s(p+2r+s)}{2} p^2. \; \] Notice
that there are two choices of $(r,s)$ which yield a simple module of
smallest possible dimension $p^2$: $((p-3)/2,1)$ and $((p-3)/2,2)$.
These weights (which parametrize a single module) lie in the second
restricted alcove, which suggests that we might instead view that
alcove as ``lowest'' and use a weight there to rewrite Jantzen's
formulas.  The alcove in question is labelled $A$ in Figure 1.

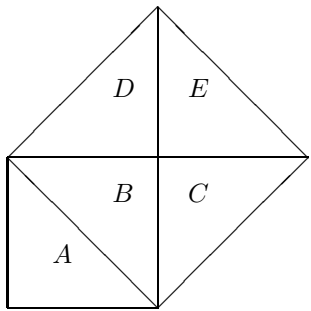
\begin{figure}[htb]\label{fig.one}
\centering

\setlength{\unitlength}{2mm}
\begin{picture}(25,25) 

\put(10,0){\line(1,0){10}}
\put(10,10){\line(1,0){20}}
\put(10,0){\line(0,1){10}}
\put(20,0){\line(0,1){20}}
\put(20,0){\line(-1,1){10}}
\put(30,10){\line(-1,1){10}}
\put(10,10){\line(1,1){10}}
\put(20,0){\line(1,1){10}}
\put(13,3){$A$}
\put(17,7){$B$}
\put(17,14){$D$}
\put(22,7){$C$}
\put(22,14){$E$}

\end{picture}

\caption{Some alcoves for type $B_2$} \end{figure} 

The four dimensions above correspond to linked weights in the
respective alcoves $A,B,C,D$.  The simple module of dimension $p^2$
corresponds to two weights in the lower left corner of alcove $A$, as
close to both vertical and horizontal walls as possible.  Moreover,
the dimension formula $s(p-2r-s)p^2/2$ for alcove $A$ corresponds in a
transparent way to the defining equations $s=0$ and $2r+s=p$ of these
two hyperplanes (if we retain Jantzen's standard coordinates).  The
two special weights minimize the numerator in this dimension formula,
giving $2p^2$; the denominator is then needed to cancel the 2.  In
this way, we can begin to imitate the interpretation of Weyl's formula
in Section \ref{sec.res}.

\subsection{} To develop further the analogy with the restricted case,
we have to rewrite the dimension formula for alcove $A$ in terms of
new ($\rho$-shifted) coordinates $(r,s)$ of a weight in this alcove.
Set \[\delta(r,s) :=s(2r+s-p)p^2/2.\] Now the idea is to apply this
formula to the ($\rho$-shifted) coordinates $(r,s)$ of an arbitrary
weight, in the spirit of Weyl's formula.  In particular the formula
yields 0 when applied to a weight in the indicated orthogonal
hyperplanes bounding alcove $A$ below.

Write briefly $\delta_A, \delta_B, \dots$ for the formal dimensions
obtained by applying $\delta$ to linked weights in alcoves $A, B,
\dots$ Denoting the corresponding simple $\uxg$-modules by $L_A, L_B,
\dots,$ we find the following pattern: \begin{eqnarray*} \dim L_A &=&
\delta_A \\ \dim L_B &=& \delta_B -\delta_A \\ \dim L_C &=& \delta_C -
\delta_B + \delta_A \\ \dim L_D &=& \delta_D - \delta_B + \delta_A
\end{eqnarray*} This in turn raises the question of the possible
existence of modules $V_\chi(\lam)$ (analogous to Weyl modules in the
case $\chi=0$) having dimensions given by the function $\delta$.
These should exist for weights lying in an appropriate collection of
alcoves (here infinite) and should be modules for $\uxg$ as well as
for $C_G(\chi)$.

The alternating sum formulas above are certainly suggestive of a
general pattern, though we usually must expect (as for $\chi=0$)
coefficients of absolute value $>1$ coming from Kazhdan--Lusztig
theory.  In our example, alcove $E$ should carry the simple module
$L_A$, but the most likely alternating sum formula will produce a
multiple of its dimension such as $3\delta_A$.  This suggests
associating to a weight in alcove $E$ a $\uxg$-module together with a
nontrivial representation of an $SL_2$-type subgroup of $C_G(\chi)$
(whose trivial representation would occur for alcove $A$).  Such a
pairing, somewhat analogous to the Springer Correspondence, would be
compatible with Lusztig's cell conjectures \cite[\S 10]{lu89}.

In any case, the main thrust of our formulation is the derivation of
diverse-looking dimension formulas from a single formula based on a
special choice of affine hyperplanes.  This much can be conjectured in
general, but the explanation for such regularity remains speculative.

\section{Kazhdan--Lusztig cells and hyperplanes}{\label{sec.cell}}

\subsection{} How can one identify suitable affine hyperplanes which
might support a version of Weyl's dimension formula for arbitrary
nilpotent $\chi$, in the spirit of the above discussion of the minimal
nilpotent orbit for type $B_2$?  An answer is suggested by Lusztig's
bijection \cite{lu89} between nilpotent orbits (in good
characteristic) and two-sided cells in the affine Weyl group (for the
Langlands dual of $G$).  As shown by Lusztig and Xi \cite{lx88}, each
two-sided cell in turn meets the dominant Weyl chamber in a
``canonical'' left cell.  In characteristic $p$ we identify the affine
Weyl group in question with $W_p$, allowing us to view the cells as
unions of $p$-alcoves.

For example, the minimal nilpotent orbit of $B_2$ corresponds to the
canonical left cell whose lower portion (beginning with alcoves $A,B,
\dots$) is the strip along one wall pictured in Figure 2.
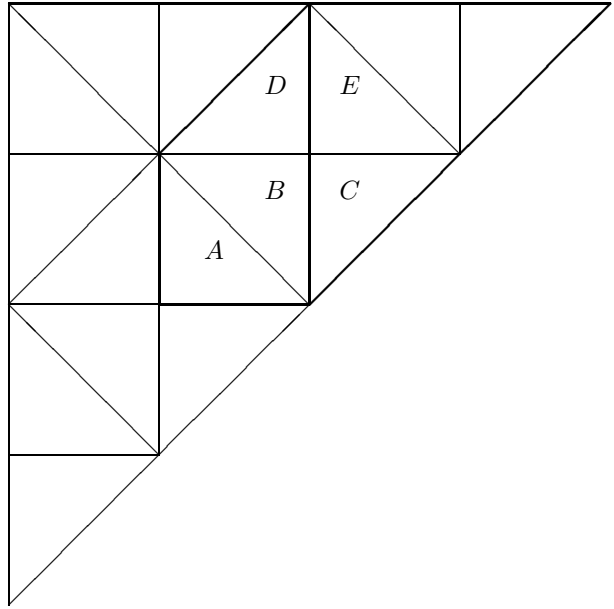
\begin{figure} \label{fig.two} \centering

\setlength{\unitlength}{2mm}
\begin{picture}(35,35)

\put(0,0){\line(1,1){20}}
\put(0,20){\line(1,1){10}}

\put(10,10){\line(-1,1){10}}
\put(20,20){\line(-1,1){20}}
\put(30,30){\line(-1,1){10}}

\put(0,0){\line(0,1){40}}
\put(10,10){\line(0,1){10}}
\put(10,30){\line(0,1){10}}
\put(20,20){\line(0,1){20}}
\put(30,30){\line(0,1){10}}

\put(0,10){\line(1,0){10}}
\put(0,20){\line(1,0){10}}
\put(0,30){\line(1,0){30}}
\put(0,40){\line(1,0){40}}

\put(13,23){$A$}
\put(17,27){$B$}
\put(22,27){$C$}
\put(17,34){$D$}
\put(22,34){$E$}

\thicklines

\put(10,20){\line(1,0){10}}
\put(10,20){\line(0,1){10}}
\put(20,20){\line(1,1){20}}
\put(10,30){\line(1,1){10}}

\end{picture}

\caption{Lower part of a canonical left cell for type $B_2$}
\end{figure}
To rewrite the previous discussion of dimensions in terms of
weights lying in these translated alcoves, we just have to
redefine the function $\delta$ by
\[\delta(r,s):=(s-p)(2r+s-2p)p^2/2.\] 
This leads to the same dimension formulas as in Sec. \ref{sec.B2}.

Empirical study of Jantzen's formulas in a variety of cases shows a
strong correlation with the hyperplanes bounding below the canonical
left cell for $\chi$ (but with hyperplanes bounding the dominant
region omitted).  This is the rationale for our reformulation of his
$B_2$ results above.

In the $B_2$ example, there are two orthogonal hyperplanes,
corresponding to an $A_1 \times A_1$ root system.  In other cases
studied one gets more complicated root systems, taking in each case
the natural hyperplanes corresponding to the associated positive roots
as a framework for the basic dimension formula.  Of course, when
$\chi=0$ we are just reverting to Weyl's formula in this way.  Such
hyperplane systems must come from various proper subsets of the
extended Dynkin diagram.

\subsection{} This type of interpretation agrees well with the location
of weights which yield $\dim \lxl = p^d$, as in the $B_2$ example.
However, this small example is oversimplified in some respects.  There
are several complicating factors in the attempt to correlate
$\uxg$-modules with cells:

\begin{itemize}

\item[(a)] It is not easy to describe geometrically the lower boundary
behavior of canonical left cells, though the work of Shi \cite{sh86}
in type $A$ (and further work by him and his associates in other
cases) provides a lot of combinatorial data.  It is clear that one
cannot in general expect to find a unique lowest alcove in a canonical
left cell.  For example, the minimal orbit for type $A_3$ (where $N=6$
and $d=3$) yields a cell with two symmetrically placed configurations
of lower hyperplanes of $A_2$ type. Here one expects three factors
corresponding to positive roots of $A_2$ in the conjectural dimension
formulas.

\item[(b)] One cannot always point to an obvious hyperplane
configuration of the right size.  An extreme case to keep in mind is
the minimal nilpotent orbit of $E_8$, where $N=120$ and $d=29$.  The
91 expected factors in a dimension formula might well arise from a
combination of the 28 positive roots in an $A_7$ root system and
another 63 positive roots in an $E_7$ root system (both found in the
extended Dynkin diagram).  It is unclear how to predict such patterns
in general, though they may be related to a duality for nilpotent
orbits studied by Sommers.  Note that for type $A$ one has a simple
version of duality (based on transpose partitions) which might suggest
a natural choice of hyperplanes.

\item[(c)] When the component group $C_G(\chi)/C_G(\chi)^\circ$ is
nontrivial, it may permute a number of nonisomorphic simple modules
having the same dimension \cite{ja99b}.  This already shows up in
subregular cases for $B_2$ or $G_2$, in a way that looks consistent
with Lusztig's conjectures in \cite[\S 10]{lu89}: each intersection of
a left cell with its inverse should correspond to an orbit of the
component group in the set of simple modules belonging to a typical
block of $\uxg$.

\end{itemize}

\subsection{} How can one construct modules $V_\chi(\lam)$ for $\uxg$
and $C_G(\chi)$ which carry dimension formulas of the shape we have
described?  One approach has been initiated by Mirkovi\'c and
Rumynin \cite{mr01}, but many technical problems remain.  The natural
starting point is the Springer fiber associated with $\chi$, whose
dimension is $N-d$.  

Ultimately all of this should connect naturally with the ideas of
Lusztig \cite{lu97,lu98,lu99a,lu99b,lu01} involving Springer fibers,
equivariant $K$-theory, affine Hecke algebras, cells, etc.
Significant progress has recently been made by Bezrukavnikov,
Mirkovi\'c, and Rumynin \cite{bmr}.

\medskip

\noindent 
{\footnotesize \emph{Acknowledgments.}
Conversations and correspondence with Jens Carsten
Jantzen have been extremely useful in formulating the ideas here,
though he should not be held responsible for my speculative
suggestions.  I am also grateful to Roman Bezrukavnikov, Paul
Gunnells, Ivan Mirkovi\'c, Jian Yi Shi, and Eric Sommers for useful
consultations.}

\bibliographystyle{amsplain}

\makeatletter \renewcommand{\@biblabel}[1]{\hfill#1.}\makeatother

\end{document}